\documentclass[12pt]{article}
\usepackage{amsmath}
\usepackage{amssymb}
\usepackage{amsmath,amssymb,amsthm,amsfonts}
{\theoremstyle{definition}

}

\newcommand{\il}[2]{\int\limits_{#1}^{#2}}

\newcommand{\ph}{\phantom{a}}
\newcommand{\phh}{\phantom{aaa}}

\newcommand{\liml}{\lim\limits_}
\newcommand{\pinf}{+\infty}

\newcommand{\no}{\eqno}
\newcommand{\bl} {\biggl\{}
\newcommand{\br} {\biggr\}}
\newcommand{\ta}{\tau}
\newcommand{\noin}{\noindent}
\newcommand{\vsk} {\vskip 10pt}
\newcommand{\al}{\alpha}
\newcommand{\be}{\beta}
\newcommand{\ga}{\gamma}

\newcommand{\et}{\eta}

\newcommand{\la}{\lambda}

\newcommand{\fr}{\frac}
\newcommand{\sist}[2]{\left\{
\begin{array}{l}
{#1}\\
\ph\\
{#2}
\end{array}
\right.}

\begin{document}

MSC 34D20

\vskip 20pt

\centerline{\bf Comparison criteria for the Abel equation of 1st kind}

\vskip 10 pt

\centerline{\bf G. A. Grigorian}

\vskip 10 pt

\centerline{0019 Armenia c. Yerevan, str. M. Bagramian 24/5}
\centerline{Institute of Mathematics of NAS of Armenia}
\centerline{E - mail: mathphys2@instmath.sci.am, \ph phone: 098 62 03 05, \ph 010 35 48 61}

\vskip 20 pt

\noindent
Abstract. Three comparison criteria for the Abel equation  of 1es kind are  proved. The  results obtained are used to obtain global solvability criteria and some criteria of existence of closed solutions  for the mentioned equation. The  results obtained  are demonstrated by examples.

\vskip 20 pt


\vskip 20 pt

{\bf  1. Introduction}.  Let $a(t), \ph b(t), \ph c(t)$ and $d(t)$ be real-valued continuous functions on $[t_0,\tau_0) \ph (t_0 < \tau_0 \le \pinf)$. Consider the Abel equation
$$
y' + a(t) y^3 + b(t) y^2 + c(t) y + d(t) = 0, \phh t_0\le t \le \tau_0. \no (1.1)
$$
This equation has important applications in questions of natural sciences (e.g., in the Friedman model, in the study of conditions for existence of limit cycles for dynamic systems [Hilbert's 16th problem]) and many works are devoted to it (see [1--5, 8--11] and cited works therein).

In this paper we prove some comparison criteria for Eq. (1.1). These criteria  we use to obtain global solvability criteria  and criteria of existence of  a closed solution for Eq. (1.1). The obtained results we demonstrate  by examples.

\vskip 20 pt

{\bf  2. Auxiliary propositions}.   Let $a_1(t), \ph b_1(t), \ph c_1(t)$ and $d_1(t)$ be real-valued continuous functions on $[t_0,\tau_0)$. Along with Eq. (1.1) consider the equation
$$
y' + a_1(t) y^3 + b_1(t) y^2 + c_1(t) y + d_1(t) = 0, \phh t_0\le t \le \tau_0. \no (2.1)
$$
Let $y(t)$ and $y_1(t)$ be solutions of Eq. (1.1) and (2.1) respectively on $[t_0,\tau_0)$. Then
$$
[y(t) - y_1(t)]' + A(t) [y(t) - y_1(t)] + B(t) = 0,
$$
where $A(t)\equiv a(t)[y^2(t) + y(t)y_1(t) + y_1^2(t)] + b(t)[y(t) + y_1(t)] + c(t), \ph B(t) \equiv [a(t) - a_1(t)] y_1^3(t) + [b(t) - b_1(t)] y_1^2(t) + [c(t) - c_1(t)]y_1(t) + d(t) - d_1(t), \ph t_0\le t < \ta_0$. We can interpret $y(t) - y_1(t)$ as a solution of the linear equation
$$
x' + A(t) x + B(t) = 0, \ph t_0 \le t <\ta_0.
$$
Then by the Cauchy formula we have
$$
y(t) - y_1(t) = \exp\bl-\il{t_0}{t}A(\ta)d\ta\br\biggl[y(t_0) - y_1(t_0) - \il{t_0}{t}\exp\bl\il{t_0}{\ta}A(s) d s\br B(\ta)d\ta\biggr], \no (2.2)
$$
$t_0 \le t < \ta_0$. Let $F(t,Y)$ be a continuous in $t$ and continuously differentiable in $Y$ vector function on $[t_0,\pinf)\times \mathbb{R}^m$. Consider the nonlinear system
$$
Y' = F(t,Y), \phh t \ge t_0. \no (2.3)
$$
Every solution $Y(t) = Y(t,t_0,Y_0)$ of this system exists either only a finite interval $[t_0,T)$ or is continuable on $[t_0,\pinf)$

\vsk

{\bf Lemma 2.1([6, p. 204, Lemma]).} {\it If a solution $Y(t)$ of the system (2.3) exists only on a finite interval $[t_0,T)$, then
$$
||Y(t)|| \to \pinf \ph \mbox{as} \ph t \to T- 0,
$$
where $||Y(t)||$ is any euclidian norm of $Y(t)$ for every fixed $t \in [t_0,T)$.
}

\phantom{aaaaaaaaaaaaaaaaaaaaaaaaaaaaaaaaaaaaaaaaaaaaaaaaaaaaaaaaaaaaaaaaaaaaaa} $\Box$

Consider the differential inequality
$$
\eta' + a(t)\eta^3 + b(t) \eta^2 + c(t) \eta + d(t) \ge 0, \phh t_0 \le t <\ta_0. \no (2.4)
$$

\vsk

{\bf Lemma 2.2.} {\it Let $\eta(t)$ be a solution of the inequality (2.4) on $[t_0,\ta_0)$ and let $y(t)$ be a solution of Eq. (1.1) on $[t_0,\ta_0)$ such that $y(t_0) \le \eta(t_0)$. Then
$$
y(t) \le \eta(t), \phh t_0 \le t < \ta_0. \no (2.5)
$$
Furthermore, if $y(t_0) < \eta(t_0)$, then
$$
y(t) < \eta(t), \phh t_0 \le t < \ta_0. \no (2.6)
$$
}

Proof. We set $\widetilde{d}(t) \equiv - \eta'(t) - a(t)\eta^3(t) - b(t)\eta^2(t) - c(t) \eta(t), \ph t_0 \le t <\ta)_0)$. Then by (2.4) we have
$$
\widetilde{\eta}(t) \le \eta(t), \phh t_0 \le t <\ta_0. \no (2.7)
$$
Obviously $\eta(t)$ is a solution of the equation
$$
y' + a(t) y^3 + b(t) y^2 + c(t) y + \widetilde{d}(t) = 0, \phh t_0\le t \le \tau_0
$$
on $[t_0,\ta_0)$. Therefore, in virtue of (2.2) we obtain
$$
y(t) - \eta(t) = \exp\bl-\il{t_0}{t} A_\eta(\ta) d \ta\br\biggl[y(t_0) - \eta(t_0) + \il{t_0}{t}\exp\bl\il{t_0}{\ta}A_\eta(s)  d s \br\biggl(\widetilde{d}(\ta) - d(\ta)\biggr) d \ta\biggr],
$$
$ t_0 \le r < \ta_0,$ where $A_\eta(t)\equiv a(t)[y^2(t) + y(t)\eta(t) + \eta^2(t)] + b(t)[y(t) + \eta(t)] + c(t), \ph t_0 \;e t < \ta_0.$ This together with (2.7) and the inequality $y(t_0) \le \eta(t_0)$ (and the inequality $y(t_0) < \eta(t_0)$) implies (2.5) (implies (2.6)). The lemma is proved.

By analogy can be proved
\vsk

{\bf Lemma 2.3.} {\it Let $\eta(t)$ be a solution of the inequality
$$
\eta' + a(t)\eta^3 + b(t) \eta^2 + c(t) \eta + d(t) \le 0, \phh t_0 \le t <\ta_0. \no (2.8)
$$
on $[t_0,\ta_0)$ and let $y(t)$ be a solution of Eq. (1.1) on $[t_0,\ta_0)$ such that $y(t_0) \ge \eta(t_0)$. Then
$$
y(t) \ge \eta(t), \phh t_0 \le t < \ta_0.
$$
Furthermore, if $y(t_0) > \eta(t_0)$, then
$$
y(t) > \eta(t), \phh t_0 \le t < \ta_0.
$$
}

\phantom{aaaaaaaaaaaaaaaaaaaaaaaaaaaaaaaaaaaaaaaaaaaaaaaaaaaaaaaaaaaaaaaaaaaaaa} $\Box$

\vsk

Consider now some cases in which Eq. (1.1) and the inequalities (2.4) and (2.8) have solutions on $[t_0,\ta_0)$.

\noin
$(I) \ph a(t)\eta^3 + b(t)\eta^2 + c(t)\eta + d(t) = a(t)(\eta - \alpha)(\eta^2 +\beta(t) \eta + \gamma(t)), \ph t_0\le t < \ta_0,$ where $\alpha(t)$ is a real-valued  continuously differentiable function on $[t_0,\ta_0), \ph \beta(t)$ and $\gamma(t)$ are real-valued continuous functions on $[t_0,\ta_0)$. Then if $\alpha'(t) \ge 0, \ph (\le 0), \ph t_0\le t <\ta_0$, the function $\alpha(t)$ is a solution of the inequality (2.4) (of the inequality (2.8)).

\noin
$(II) \ph a(t)\eta^3 + b(t)\eta^2 + c(t)\eta + d(t) = a(t)(\eta - \alpha)(\eta-\beta(t))(\eta - \gamma(t)), \ph a(t) \ne 0, \ph t_0\le t < \ta_0,$ where $\alpha(t), \ph \beta(t)$ and $\gamma(t)$ are some real-valued continuous functions on $[t_0,\ta_0), \ph \alpha(t) \le \beta(t) \le \gamma(t), \ph t_0 \le t <\ta_0$.

\noin
$(II_1)$ \ph Assume $a(t) > 0, \ph t_0 \le t <\ta_0$ and $\alpha(t) \le \eta_0 \le \beta(t), \ph (\beta(t) \le \eta_1 \le \gamma(t)), \ph t_0\le t <~ \ta_0$ for some $\eta_0 = const, \ph \eta_1 = const$. Then $\eta(t)\equiv \eta_0 \ph (\eta(t) \equiv \eta_1)$ will be  a solution of the inequality (2.4) (of the inequality (2.8)).

\noin
$(II_2)$ \ph Assume $a(t) < 0, \ph t_0 \le t <\ta_0$ and $\alpha(t) \le \eta_0 \le \beta(t), \ph (\beta(t) \le \eta_1 \le \gamma(t)), \ph t_0\le t < ~\ta_0$ for some $\eta_0 = const, \ph \eta_1 = const$. Then $\eta(t)\equiv \eta_0 \ph (\eta(t) \equiv \eta_1)$ will be  a solution of the inequality (2.8) (of the inequality (2.4)).

\vsk

{\bf Theorem 2.1.} {\it Let the following conditions be satisfied.

\noin
(III) \ph $a(t) > 0$ almost everywhere on $[t_0,\ta_0)$,

\noin
(IV) \ph $\frac{b^2(t)}{a(t)}$ is locally integrable on $[t_0,\ta_0)$.

\noin
Then for every $\ga \in\mathbb{R}$ the solution $y(t)$ of Eq. (1.1) with $y(t_0) =\ga$ exists on $[t_0,\ta)$.
}

Proof. Let $y(t)$ be a solution of Eq. (1.1) with $y(t_0) = \ga \in \mathbb{R}$. and let $[t_0,t_1)$ be its maximum existence interval. We must show that
$$
t_1 = \ta_0. \no (2.9)
$$
By (1.1) we have
$$
y'(t) + A_0(t) y(t) + d(t) = 0, \phh t \in [t_0,t_1),
$$
where $A_0(t) \equiv a(t) y^2(t) + b(t) y(t) + c(t), \ph t \in [t_0,t_1)$. Then by the Cauchy formula we have
$$
y(t) =\ga \exp\bl-\il{t_0}{t}A_0(\ta) d \ta\br- \il{t_0}{t}\exp\bl-\il{\ta}{t} A_0(s) d s\br d(\ta) d\ta, \phh t \in [t_0,t_1). \no (2.10)
$$
It follows from the condition (III) that $A_0(t) \ge c(t) - \fr{b^2(t)}{4a(t)}$ almost everywhere on $[t_0,t_1)$. Then from (2.10) we get
$$
|y(t)| \le |\ga| \exp\bl -\il{t_0}{t}\Bigl[c(\ta) -\frac{b^2(\ta)}{4a(\ta)}\Bigr]d\ta\br + \il{t_0}{t}\exp\bl-\il{\ta}{t} \Bigl[c(s) - \fr{b^2(s)}{4a(s)}\Bigr] d s\br |d(\ta)| d\ta, \no (2.11)
$$
$t \in [t_0,t_1).$
This inequality with the condition (IV) implies that $y(t)$ is bounded on $[t_0,t_1)$. Then by Lemma 2.1 the supposition $t_1 < \ta_0$ leads to a contradiction. Hence, (2.9) holds. The theorem is proved.

{\bf Remark 2.1.} {\it In the case $\ta_0 < \pinf$ Theorem 2.1 remains valid if $[t_0,\ta_0)$ is replaced  in it by $[t_0,\ta_0]$. Indeed,  it follows from the  integrability of the function $\fr{b^2(t)}{a(t)}$ on $[t_0,\ta_0]$ and from (2.11) with $t_1 = \ta_0$ that $y(t)$ is bounded on $[t_0,\ta_0)$. Then  the equality
$$
y(t) = y(t_0) - \il{t_0}{t} [a(\ta) y^3(\ta) + b(\ta) y^2(\ta) + c(\ta) y(\ta) + d(\ta)] d \ta, \ph t \in [t_0,\ta),
$$
 which holds by (1.1), implies that $y(t)$ is continuable on $[t_0,\ta_0]$ as a solution of Eq. (1.1).
}

\vsk

{\bf 3. Comparison criteria.}

\vsk

{\bf Definition 3.1.} {\it An interval $[t_0,t_1) \ph (t_0 < t_1 \le \ta_0)$ is called the maximum existence interval for a solution $y(t)$ of Eq. (1.1) if $y(t)$ exists on $[t_0,t_1)$ and cannot be continued to the right from $t_1$ as a solution of Eq. (1.1).
}

\vsk

{\bf Theorem 3.1.} {\it Let $y_1(t)$ be a solution of Eq. (2.1) on $[t_0,\ta_0)$ and let $\et(t)$ be a solution of the inequality (2.4) on $[t_0,\ta_0)$ with $\et(t_0) \ge y_1(t_0)$. Moreover, 
let the following conditions be satisfied

\noin
(I) \ph  $a(t) < 0$ almost everywhere on $[t_0,\ta_0),$ and the function $\fr{b^2(t)}{a(t)}$ is locally integrable on $[t_0,\ta_0)$

\noin
(II) \ph $\ga - y_1(t_0) + \il{t_0}{t}\exp\bl\il{t_0}{\ta}\Bigl[c(s) - \frac{b^2(s)}{a(s)}\Bigr]ds\br\biggl[(a_1(\ta) - a(\ta))y_1^3(\ta) + (b_1(\ta) - b(\ta))y_1^2(\ta) +\linebreak  + (c_1(\ta) - c(\ta)) y_1(\ta) + d_1(\ta) - d(\ta)\biggr] d \ta \ge 0, \ph t_0 \le t < \ta_0$ for some $\ga \in [y_1(t_0),\et(t_0)]$.

\noin
Then every solution $y(t)$ of Eq. (1.1) with $y(t_0) \in [\ga,\et(t_0)]$ exists on $[t_0,\ta_0)$ and
$$
y_1(t) \le y(t) \le \et(t), \phh t_0\le t < \ta_0. \no (3.1)
$$
Furthermore, if $y_1(t_0) < y(t_0) \ph (y(t_0) < \et(t_0))$ then
$$
y_1(t) < y(t) \phh (y(t) < \et(t)), \phh t_0\le t < \ta_0. \no (3.2)
$$
}

Proof. Let $y(t)$ be a solution of Eq. (1.1) with $y(t_0) \in [\ga,\et(t_0)]$ and let $[t_0,t_1)$ be its maximum existence interval. Then by virtue of Lemma 2.2
$$
y(t) \le \et(t), \phh t_0 \le t < t_1 \no (3.3)
$$
and if $y(t_0) < \et(t_0),$ then
$$
y(t) < \et(t), \phh t_0 \le t < t_1. \no (3.4)
$$
In virtue of (2.2) we have
$$
y(t) - y_1(t) = \exp\bl-\il{t_0}{t}A_1(\ta)d \ta\br\biggl[y(t_0) - y_1(t_0) + \il{t_0}{t}\exp\bl\il{t_0}{\ta} A_1(s) d s\br B_1(\ta) d \ta\biggr], \no (3.5)
$$
$t_0 \le t < t_1$, where $A_1(t) \equiv a(t)[y^2(t) + y(t) y_1(t) + y_1^2(t)] + b(t)[y(t) - y_1(t)] + c(t), \ph B_1(t) \equiv [a_1(t) - a(t)]y_1^3(t) + [b_1(t) - b(t)] y_1^2(t) + [c_1(t) - c(t)] y_1(t) + d_1(t) - d(t), \ph t_0 \le t < t_1.$ It is not difficult to derive from the condition (I) that
$$
A_1(t) \le c(t) - \fr{b^2(t)}{a(t)}, \ph \mbox{amost everywhere on} \ph [t_0,t_1). \no (3.6)
$$
We have
$$
\il{t_0}{t}\exp\bl\il{t_0}{\ta}A_1(s) d s\br B_1(\ta) d \ta = \il{t_0}{t} G(\ta)exp\bl\il{t_0}{\ta}\Bigl[c(s) - \fr{b^2(s)}{a(s)}\Bigr]ds\br B_1(\ta) d \ta,  \no (3.7)
$$
$t_0\le t < t_1$, where $G(t)\equiv\exp\bl\il{t_0}{t}[A_1(s) + \fr{b^2(s)}{a(s)} - c(s)]d s\br, \ph t_0 \le t < t_1.$
It follows from (3.6) that $G(t)$ is a monotonically non increasing function on $[t_0,t_1)$. Then by the mean value theorem for integrals (see [7, p. 869]) from (3.7) we obtain
$$
\il{t_0}{t}\exp\bl\il{t_0}{\ta}A_1(s) d s\br B_1(\ta) d \ta = \il{t_0}{\al(t)}\exp\bl\il{t_0}{\ta}\Bigl[c(s) - \fr{b^2(s)}{a(s)}\Bigr] d s\br B_1(\ta) d \ta, \ph t_0\le t \le t_1,
$$
For some $\al(t) \in [t_0,t).$ This together with the condition (II) and (3.5) implies that
$$
y_1(t) \le y(t), \phh t_0 \le t < t_1, \no (3.8)
$$
and if $y_1(t_0) < y(t_0)$ then
$$
y_1(t) <  y(t), \phh t_0 \le t < t_1.
$$
This relation with (3.3), (3.4) and (3.8) shows that the proof of the theorem will be completed if we show that
$$
t_1 = \ta_0. \no (3.9)
$$
Suppose $t_1 < \ta_0$. Then it follows from (3.3) and (3.8) that $y(t)$ is bounded on $[t_0,t_1)$. Then by virtue of Lemma 2.2 $[t_0,t_1)$ is not the maximum existence interval for $y(t)$. We have obtained a contradiction, completing the proof of the theorem.

Using Lemma 2.3 instead of Lemma 2.2 by analogy with the proof of Theorem 3.1 can be proved the following theorem.

\vsk

{\bf Theorem 3.2.}. {\it Let $y_1(t)$ be a solution of Eq. (2.1) on $[t_0,\ta_0)$ and let $\et(t)$ be a solution of the inequality (2.8) on $[t_0,\ta_0)$ with $\et(t_0) \le y_1(t_0)$. Moreover, let 
 the  condition (I) of  Theorem 3.1
and the condition

\noin
$\ga - y_1(t_0) + \il{t_0}{t}\exp\bl\il{t_0}{\ta}\Bigl[c(s) - \frac{b^2(s)}{a(s)}\Bigr]d s\br\biggl[(a_1(\ta) - a(\ta))y_1^3(\ta) + (b_1(\ta) - b(\ta))y_1^2(\ta) +\linebreak  + (c_1(\ta) - c(\ta)) y_1(\ta) + d_1(\ta) - d(\ta)\biggr] d \ta \le 0, \ph t_0 \le t < \ta_0$ for some $\ga \in [\et(t_0),y_1(t_0)]$

\noin
be satisfied.

\noin
Then every solution $y(t)$ of Eq. (1.1) with $y(t_0) \in [\et(t_0),\ga]$ exists on $[t_0,\ta_0)$ and
$$
y_1(t) \le y(t) \le \et(t), \phh t_0\le t < \ta_0.
$$
Furthermore, if $\et(t_0) < y(t_0) \ph (y(t_0) < y_1(t_0))$ then
$$
\et(t) < y(t) \ph (y(t) < y_1(t)), \ph t_0 \le t < \ta_0.
$$
}

\phantom{aaaaaaaaaaaaaaaaaaaaaaaaaaaaaaaaaaaaaaaaaaaaaaaaaaaaaaaaaaaaaaaaaaaaaa} $\Box$

Let $a_2(t), \ph b_2(t), \ph c_2(t)$ and $d_2(t)$ be real-valued continuous functions on $[t_0,\ta_0)$. Consider the equation
$$
y' + a_2(t) y^3 + b_2(t) y^2 + c_2(t) y + d_2(t) = 0, \phh t_0 \le t < \ta_0. \no (3.10)
$$

\vsk

{\bf Theorem 3.3.} {\it Let $y_1(t)$ and $y_2(t)$ be solutions of Eq. (2.1) and Eq. (3.10) respectively on $[t_0,\ta_0)$ such that $y_1(t_0) \le y_2(t_0)$, 
 and let the following conditions be satisfied

\noin
(I) \ph $a(t) < 0$ almost everywhere on $[t_0,\ta_0)$ and the function $\fr{b^2(t)}{a(t)}$ is locally integrable on $[t_0,\ta_0)$,

\noin
(II) \ph $\ga_1 - y_1(t_0) + \il{t_0}{t}\exp\bl\il{t_0}{\ta}\Bigl[c(s) - \frac{b^2(s)}{a(s)}\Bigr]d s\br\biggl[(a_1(\ta) - a(\ta))y_1^3(\ta) + (b_1(\ta) - b(\ta))y_1^2(\ta) +\linebreak  + (c_1(\ta) - c(\ta)) y_1(\ta) + d_1(\ta) - d(\ta)\biggr] d \ta \ge 0, \ph t_0 \le t < \ta_0$ for some $\ga_1 \in [y_1(t_0),y_2(t_0)]$,

\noin
(III) \ph $\ga_2 - y_2(t_0) + \il{t_0}{t}\exp\bl\il{t_0}{\ta}\Bigl[c(s) - \frac{b^2(s)}{a(s)}\Bigr]d s\br\biggl[(a_2(\ta) - a(\ta))y_2^3(\ta) + (b_2(\ta) - b(\ta))y_2^2(\ta) +\linebreak  + (c_2(\ta) - c(\ta)) y_2(\ta) + d_2(\ta) - d(\ta)\biggr] d \ta \le 0, \ph t_0 \le t < \ta_0$ for some $\ga_2 \in [\ga_1,y_2(t_0)].$

\noin
Then every solution $y(t)$ of Eq. (1.1) with $y(t_0) \in [\ga_1,\ga_2]$  exists on $[t_0,\ta_0)$ and
$$
y_1(t) \le y(t) \le y_2(t), \phh t_0 \le t < \ta_0. \no (3.11)
$$
Furthermore, if $y_1(t_0) < y(t_0) \ph (y(t_0) < y_2(t_0))$ then
$$
y_1(t) < y(t) \ph (y(t) < y_2(t)), \phh t_0 \le t < \ta_0. \no (3.12)
$$
}

Proof.  Let $y(t)$ be a solution of Eq. (1.1) with $y(t_0) \in [\ga_1,\ga_2]$ and let $[t_0,t_1)$ be its maximum existence interval. By (2.2) we have
$$
y(t) - y_k(t) = \exp\bl-\il{t_0}{t}A_k(\ta)d\ta\br\biggl[y(t_0) - y_k(t_0) - \il{t_0}{t}\exp\bl\il{t_0}{\ta}A_k(s) d s\br B_k(\ta)d\ta\biggr], \no (3.13)
$$
where
$A_k(t)\equiv a(t)[y^2(t) + y(t)y_k(t) + y_k^2(t)] + b(t)[y(t) + y_k(t)] + c(t), \ph B(t) \equiv [a(t) - a_k(t)] y_k^3(t) + [b(t) - b_k(t)] y_k^2(t) + [c(t) - c_k(t)]y_k(t) + d(t) - d_k(t), \ph t_0\le t < t_1, \ph k=1.2.$. The condition (I) implies that
$$
A_k(t) - c(t) + \fr{b^2(t)}{a(t)} \le 0, \ph \mbox{almost everywhere on} \ph [t_0,t_1).
$$
Then by mean value theorem for integral the following equality holds
$$
\il{t_0}{t}\exp\bl\il{t_0}{\ta}A_k(s) d s\br B_k(\ta) d \ta = \il{t_0}{\al_k(t)}\exp\bl\il{t_0}{\ta}\Bigl[a(s) - \fr{b^2(s)}{a(s)}\Bigr]d s\br B_k(\ta) d \ta, \ph t_0 \le \al_k(t) \le t < t_1,
$$
$k=1,2.$ This together with (3.13) and the conditions (II) and (III) implies that
$$
y_1(t) \le y(t) \le y_2(t), \phh t_0 \le t < t_1, \no (3.14)
$$
and if $y_1(t_0) < y(t_0) \ph (y(t_0) < y_2(t_0))$ then
$$
y_1(t) < y(t) \ph (y(t) < y_2(t)), \phh t_0 \le t < t_1. \no (3.15)
$$
Show that
$$
t_1 = \ta_0. \no (3.16)
$$
Suppose $t_1 < \ta_0.$ Then from (3.14) it follows that $y(t)$ is bounded on $[t_0,t_1)$. By Lemma 2.1 it follows from here that $[t_0,t_1)$ is not the maximum existence interval for $y(t)$. We have obtained a contradiction, proving (3.18). The relations (3.14)-(3.16) imply the relations (3.11) and (3.12). The theorem is proved.

\vsk

{\bf Theorem 3.4.} {\it Let $y_1(t)$ be a solution of Eq. (2.1) on $[t_0,\ta_0)$ and let the following conditions be satisfied.

\noin
(I) \ph $a(t) > 0$ almost everywhere on $[t_0,\ta_0)$,

\noin
(IV) \ph $[a_1(t) - a(t)] y_1^3(t) + [b_1(t) - b(t)]y_1^2(t) + [c_1(t) - c(t)]y_1(t) + d_1(t) - d(t) \ge 0 \ph (\le 0), \ph t \in [t_0,\ta_0)$.

\noin
Then any solution $y(t)$ of Eq. (1.1) with $y(t_0) \le y_1(t_0) \ph (y(t_0) \ge y_1(t_0))$ satisfies the inequality.
$$
y(t) \ge y_1(t) \ph (y(t) \le y_1(t)), \phh t \in [t_0,\ta_0)
$$
and if $y(t_0) > y_1(t_0) \ph (y(t_0) < y_1(t_0))$, then $y(t) > y_1(t) \ph (y(t) < y_1(t)), \ph t \in [t_0,\ta_0)$.
}

Proof. According to Theorem 2.1 the existence of $y(t)$ on $[t_0,\ta_0)$ follows from the condition (I). Then by (2.2) we have
$$
y(t) - y(t_0) = \exp\bl-\il{t_0}{t}A_1(\ta)d \ta\br\Bigl[y(t_0) - y_1(t_0) + \il{t_0}{t}\exp\bl\il{t_0}{\ta}A_1(s) d s\br B_1(\ta) d \ta\Bigr], \ph t \in [t_0,\ta_0),
$$
where $A_1(t) \equiv a(t)[y^2(t) + y(t) y_1(t) + y_1^2(t)] + b(t)[y(t) - y_1(t)] + c(t), \ph B_1(t) \equiv [a_1(t) - a(t)]y_1^3(t) + [b_1(t) - b(t)] y_1^2(t) + [c_1(t) - c(t)] y_1(t) + d_1(t) - d(t), \ph t \in [t_0,\ta_0).$ This together with the conditions (IV) and $y(t_0) \le y_1(t_0) \ph (y(t_0) \ge y_1(t_0))$ implies $y(t) > y_1(t) \ph (y(t) < y_1(t)), \ph t \in [t_0,\ta_0)$, and if $y(t_0) > y_1(t_0) \ph (y(t_0) < y_1(t_0))$, then $y(t) > y_1(t) \ph (y(t) < y_1(t)), \ph t \in [t_0,\ta_0)$. The theorem is proved.

\vsk
{\bf Theorem 3.5.} {\it Let $y_1(t)$ and $y_2(t)$ be solutions of Eq. (2.1) and Eq. (3.10) respectively on $[t_0,\ta_0)$ such that $y_1(t_0) \le y_2(t_0)$ and let the following conditions be satisfied

\noin
(IV) \ph $(a_1(t) - a(t)) y_1^3(t) + (b_1(t) - b(t)) y_1^2(t) + (c_1(t) - c(t)) y_1(t) + d_1(t) - d(t) \ge 0, \ph t \in [t_0,\ta_0),$

\noin
(V) \ph $(a_2(t) - a(t)) y_2^3(t) + (b_2(t) - b(t)) y_2^2(t) + (c_2(t) - c(t)) y_2(t) + d_2(t) - d(t) \le 0, \ph t \in [t_0,\ta_0).$

\noin
Then every solution $y(t)$ of Eq. (1.1) with $y(t_0) \in [y_1(t_0),y_2(t_0)]$ exists on $[t_0,\ta_0)$ and
$$
y_1(t) \le y(t) \le y_2(t), \phh t \in [t_0,\ta_0).
$$
Furthermore, if $y_1(t_0) < y(t_0)\ph (y(t_0) < y_2(t_0))$ then
$$
y_1(t) < y(t) \phh (y(t) < y_2(t)) \phh t \in [t_0,\ta_0).
$$
}

Proof. Let $y(t)$ be a solution of Eq. (1.1) with $y(t_0) \in [y_1(t_0),y_2(t_0)]$ and let $[t_0,t_1)$ be its maximum existence interval. Then by (2.2) for $y(t), \ph y_1(t)$ and $y_2(t)$ the equalities (3.13) hold on $[t_0,t_1)$. It follows from (3.13) and from the conditions (IV) and (V) of the theorem that
$$
y_1(t) \le y(t) \le y_2(t), \phh t \in [r_0,t_1) \no (3.17)
$$
and if $y_1(t_0) < y(t_0)\ph (y(t_0) < y_2(t_0))$ then
$$
y_1(t) < y(t) \phh (y(t) < y_2(t)) \phh t \in [t_0,1).
$$
Therefore, to complete the proof of the theorem it remains to show that $t_1 = \ta_0$. Suppose $t_1 < \ta_0$. Then it follows from (3.17) that $y(t)$ is bounded on $[t_0,\ta_0)$. By Lemma 2.1 it follows from here that $[t_0,t_1)$ is not the maximum existence interval for $y(t)$. We have obtained a contradiction, completing the proof of the theorem.

\vsk
{\bf Remark 3.1.} {\it In the case $\ta_0<\pinf$ Theorems 3.1--3.5 remain valid if  $[t_0,\ta_0)$ is replaced in them by $[t_0,\ta_0]$.
}

\vsk

{\bf Example 3.1.} {\it Consider the equation
$$
y' - y^3 + y - \la   \sin ^2 t = 0, \phh t \ge t_0. \no (3.18)
$$
where $\la = const \in [0,\fr{2\sqrt{3}}{9}].$ One can easily show that $\et(t) \equiv \fr{\sqrt{3}}{3}$ is a solution of the differential inequality
$$
\et' - \et^3 + \et - \la \sin^2 t \ge 0, \phh t \ge 0. \no (3.19)
$$
on $[0, \pinf)$. Consider the equation
$$
y' - y^3 + y = 0, \phh t \ge 0. \no (3.20)
$$
Obviously $y_1(t) \equiv 0$ is a solution of this equation on $[0,\pinf)$. Then applying Theorem~ 3.1 to (3.18)-(3.20) we conclude that every solution $y(t)$ of Eq. (3.18) with $y(0)\in [0,\fr{\sqrt{3}}{3}]$ exists on $[0,\pinf)$ and
$$
0\le y(t) \le \fr{\sqrt{3}}{3}, \phh t \ge 0.
$$
Moreover, if $y(0) > 0 \ph (y(0) <  \fr{\sqrt{3}}{3})$ then
$$
y(t) > 0 \ph \Bigl(y(t) <  \fr{\sqrt{3}}{3}\Bigr), \phh t \ge 0.
$$
}

\vsk

{\bf Example 3.2.} {\it Using Theorem 3.2 by analogy with Example 3.1 one can show that every solution of the equation
$$
y' - y^3 + y + \la \sin^2 t = 0, \phh t \ge 0 \ph \Bigl(\la \in \Bigl[0, \fr{2\sqrt{3}}{9}\Bigr]\Bigr)
$$
with $-\fr{\sqrt{3}}{3} \le y(0) \le 0$ exists on $[0,\pinf)$ and
$$
-\fr{\sqrt{3}}{3} \le y(t) \le 0 \phh t \ge 0,
$$
and if  $y(0) >  - \fr{\sqrt{3}}{3}   \ph (y(0) < 0)$, then
$$
y(t) >  - \fr{\sqrt{3}}{3}   \ph (y(t) < 0) , \phh t \ge 0.
$$
}

\vsk
{\bf Example 3.3.} {\it Consider the equations
$$
y' - y^3 + 3 y^2 +3 y - 3 - \mu \sin t = 0, \phh t \ge 0, \ph \mu = const,  \ph |\mu| \le 2, \no (3.21)
$$

$$
y' + y^3 + 3 y^2 + 3 y  + 1  = 0, \phh t \ge 0. \no (3.22)
$$

$$
y' - y^3 + 3 y^2 - 3 y  + 1 = 0, \phh t \ge 0, \no (3.23)
$$

Obviously $y_1(t)\equiv -1, \ph t \ge 0$ is a solution of Eq. (3.22) on $[0,\pinf)$ and $y_2(t)\equiv  1, \ph t \ge 0$
is a solution of Eq. (3.23) on $[0,\pinf)$. Then using Theorem 3.3 to the equations (3.21)-(3.23) we conclude that every solution $y(t)$ of Eq. (3.21) with $|y(0)| \le 1$ exists on $[0,\pinf)$ and
$$
|y(t)| \le 1, \phh t \ge 0,
$$
and if $y(0) > -1 \ph (y(0) <1)$ then
$$
y(t) > -1 \ph (y(t) <1), \phh t \ge 0.
$$
}

\vsk
{\bf Example 3.4.} {\it Consider the equation
$$
y' + \mu \sin t y^3 + 3 y^2 + 3 y + \la(t) = 0, \phh t \ge 0,
$$
where $\mu = const, \ph |\mu| \le 1, \ph \la(t)$ is a continuous function on $[0,\pinf)$ such that $-5 \le \la(t) \le -2. \ph t \ge 0$. Then using Theorem 3.5 to this and the equations (3.22) and (3.23) one can easily show that every solution $y(t)$ of the last equation with $|y(0)| \le 1$ exists on $[0,\pinf)$ and $|y(t)| \le 1, \ph t \ge 0.$
}

\vsk

{\bf 4. Global solvability criteria}. In this section, we derive some global criteria for the existence of the Abel equation (1.1) as simple consequences of the obtained comparison criteria. In section 2 it was sown that if $a(t) >0$ almost everywhere on $[t_0,\pinf)$, then for every initial value $\la \in \mathbb{R}$ every solution $y(t)$ of Eq. (1.1) with $y(t_0) = \la$
exists on $[t_0,\pinf)$. In this section will be considered the case $a(t) < 0$ almost everywhere on $[t_0,\pinf)$.

Let $t_0< t_` < \dots$ be a finite or infinite sequence such that $t_k \in [t_0,\pinf]$. We assume that if $\{t_k\}$ is finite then then the greatest value among $t_k$ equals $\pinf$ and if  $\{t_k\}$ is infinite then $\lim\limits_{k \to \pinf} t_k = \pinf.$

\vsk
{\bf Theorem 4.1.} {\it Let 
 the following conditions be satisfied

\noin
(i) \ph $a(t) < 0$ almost everywhere on $[t_0,\pinf)$  and he function $\fr{b^2(t)}{a(t)}$ is locally integrable on $[t_0,\pinf)$,

\noin
(ii) the inequality (2.4) has a nonnegative solution $\et(t)$ on $[t_0,\pinf)$,
$$
(iii)  \phantom{aaaaaaa}  \et(t_k) - \il{t_k}{t}\exp\bl\il{t_k}{\ta}\Bigl[c(s) - \fr{b^2(s)}{a(s)}\Bigr]d s\br d(\ta) d \ta \ge 0, \ph t \in [t_k,t_{k+1}), \ph k=0,1,\dots. \phantom{aaaaaaaaaaaaaaaaaaaaaaaaaa}
$$
Then every solution $y(t)$ of Eq. (1.1) with $0\le y(t_0) \le \et(t_0)$ exists on $[t_0,\pinf)$ and
$$
0\le y(t) \le \et(t), \phh t \ge t_0.
$$
}

Proof. Let us take $a_1(t)\equiv a(t), \ph b_1(t) \equiv b(t), \ph c_1(t)\equiv c(t), \ph d_1(t)\equiv 0, \ph t \ge t_0, \ph t_0 = t_k,ph \ta_0 = t_{k+1}$. Then using Theorem 3.1. we conclude that for each $k = 0,1, \dots$ Eq. (1.1) has a solution $y_k(t)$ on $[t_k,t_{k+1})$ with with $y_k(t_k) \in [0,\et(t_k)]$ and $0\le y_k(t) \le \et(t), \ph t \in [t_k,t_{k+1}).$. Therefore $y_k(t)$ is bounded on $[t_k,t_{k+1})$. Show that $y_k(t_{k+1} - 0)$ exists for each $k =0,1,\dots$. By (1.1) we have
$$
y_k(t) = y_k(t_k) + \il{t_k}{t}[a(\ta)y_k^3(\ta) + b(\ta) y_k^2(\ta) + c(t) y_k(\ta) + d(\ta)]d\ta, \ph t \in [t_k,t_{k+1}), \ph k=0,1,\dots.
$$
It follows from here and from the boundedness of $y_k(t)$ on $[t_k,t_{k+1})$ the existence of $y_k(t_{k+1} - 0)$. Then to complete the proof it is enough take chose $y_{k+1}(t_{k+1}) = y_k(t_{k+1} - 0), \ph k = 0,1,\dots$ (as far as $y_k(t_k+1) - 0) \le \et(t_{k+1}), \ph k =0,1,\dots$ and determine $y(t)$ as a solution of Eq. (1.1) as follows
$$
y(t) = y_k(t), \phh t \in [t_k,t_{k+1}), \phh k=0,1,\dots.
$$
The theorem is proved.

By analogy using Theorem 3.2 it can be proved

\vsk

{\bf Theorem 4.2.} {\it Let   the following conditions be satisfied

\noin
 \ph $a(t) < 0$ almost everywhere on $[t_0,\pinf)$ and the function $\fr{b^2(t)}{a(t)}$ is locally integrable on $[t_0,\pinf)$,

\noin
 the inequality (2.8) has a non positive solution $\et(t)$ on $[t_0,\pinf)$,
$$
\phantom{aaaaa} \et(t_k) - \il{t_k}{t}\exp\bl\il{t_k}{\ta}\Bigl[c(s) - \fr{b^2(s)}{a(s)}\Bigr]d s\br d(\ta) d \ta \le 0, \ph t \in [t_k,t_{k+1}), \ph k=0,1,\dots. \phantom{aaaaaaaaaaaaaaaaaaaaaaaaa}
$$
Then every solution $y(t)$ of Eq. (1.1) with $\et(t_0) \le y(t_0) \le 0$ exists on $[t_0,\pinf)$ and
$$
\et(t) \le y(t) \le 0, \phh t \ge t_0.
$$
}

\phantom{aaaaaaaaaaaaaaaaaaaaaaaaaaaaaaaaaaaaaaaaaaaaaaaaaaaaaaaaaaaaaaaaaaaaaa} $\Box$

Indicate two case in which the condition (ii) of Theorem 4.1 is satisfied

\noin
1) $a(t) y^3 + b(t) y^2 + c(t) y + d(t) = a(t)(y - \al(t))(y^2 + \be(t) y + \ga(t)), \ph t \ge t_0, \ph y \in \mathbb{R}$, where $\al(t), \ph \be(t)$ and $\ga(t)$ are real-valued functions, $\al(t)$ is nonnegative nondecreasing continuously differentiable on $[t_0,\pinf)$.

\noin
2) $a(t) y^3 + b(t) y^2 + c(t) y + d(t) = a(t)(y - \al(t))(y - \be(t))(y - \ga(t)), \ph t \ge t_0$, where $\al(t)\le \be(t) \le \ga(t), \ph t \ge t_0$ are real-valued functions, $\be(t) \le \la(t) \le \ga(t), \ph t\ge t_0, \ph \la(t)$ is a nonnegative nondecreasing function on $[t_0,\pinf)$ (then $\et(t)\equiv \la(t)$ is a solution of the inequality (2.4)).

Similar cases can be indicated for  satisfiability of the condition (iv) of Theorem 4.2.
\vsk

{\bf 5. Existence of a closed solution.}

{\bf Definition 5.1 (see [11]).} {\it A solution $y(t)$ of Eq. (1.1), existing on an interval $[t_0,T]$, is called closed on $[t_0,T]$, if $y(t_0) = y(T)$,
}

The problem of finding a closed solution of Eq. (1.1) is closely related to Hilbert's 16th problem and many papers are devoted to it (see [11] and the papers cited there). For polynomial differential equations (for the generalized Abel equation), some criteria for the existence of a closed solution were obtained in [11]. In this section, based on the results of section 3, we prove some criteria for the existence of a closed solution to Eq.  (1.1). These criteria are not covered by the results of work [11].

{\bf Theorem 5.1.} {\it Let the following conditions be satisfied

\noin
(1) \ph $a(t) > 0$ almost everywhere on $[t_0,T]$,

\noin
(2) \ph $\fr{b^2(t)}{a(t)}$ is integrable on $[t_0,T]$,

\noin
(3) \ph  $\il{t_0}{T}\Bigl[c(t) - \fr{b^2(t)}{4 a(t)}\Bigr]d t > 0,$

\noin
(4) \ph $d(t) \le 0,  \ph t \in [t_0,T]$.

\noin
Then Eq. (1.1) has a nonnegative closed solution  on $[t_0,T]$.
}

Proof. Let $y_\ga(t)$ be a solution of Eq. (1.1) with $y_\ga(t_0) = \ga \ge 0$. By Theorem 3.4 and Remark 3.1 it follows from the conditions (1) and (2) that $y_\ga(t)$ exists on $[t_0,T]$ and is nonnegative. Show that there exists $\ga_0$ such that
$$
y_{\ga_0}(t_0) \ge y_{\ga_0}(T). \no (5.1)
$$
By (1.1) we have
$$
y_\ga'(t) + A_\ga(t) y_\ga(t) + d(t) = 0, \phh t \in [t_0,T],
$$
where $A_\ga(t) \equiv a(t) y_\ga^2(t) + b(t) y_\ga(t) + c(t), \ph t \in [t_0,T]$. Then by the Cauchy formula we get
$$
y_\ga(t) = \ga \exp\bl-\il{t_0}{t}A_\ga(\ta) d \ta\br - \il{t_0}{t}\exp\bl-\il{\ta}{t}A_\ga(s) d s\br d(\ta) d \ta, \phh t \in [t_0,T]. \no (5.2)
$$
It is not difficult to verify that $A_\ga(t) \ge c(t) - \fr{b^2(t)}{4 a(t)}$ almost everywhere on $[t_0,T].$ Then the conditions (2) and (4) imply
$$
y_\ga(t) \le \ga \exp\bl-\il{t_0}{t}\Bigl[c(\ta) - \fr{b^2(\ta)}{4 a(\ta)}\Bigr]d \ta \br  + \il{t_0}{t}\exp\bl-\il{\ta}{t}\Bigl[c(s) - \fr{b^2(s)}{4 a(s)}\Bigr]d s\br |d(\ta)| d \ta,  \no (5.3)
$$
$t \in [t_0,T].$ It follows from the condition (4) that
$$
1 - \exp\bl-\il{t_0}{T}\Bigl[c(\ta) - \fr{b^2(\ta)}{4 a(\ta)}\Bigr]d\ta\br > 0.
$$
This together with (5.3) implies that if
$$
\ga_0 \ge \fr{\il{t_0}{T}\exp\bl-\il{\ta}{T}\Bigl[c(s) - \fr{b^2(s)}{4 a(s)}\Bigr]d s\br|d(\ta)| d \ta}{1 - \exp\bl-\il{t_0}{T}\Bigl[c(\ta) - \fr{b^2(\ta)}{4 a(\ta)}\Bigr]d\ta\br},
$$
then (5.1) holds. Define two sequences $\{\xi_n\}_{n=1}^{\pinf}$ and  $\{\et_n\}_{n=1}^{\pinf}$ by induction on $n$ as follows:

\noin
$\xi_1 = 0, \ph \et_1 = \ga_0, \ph \xi_{n+1} =\xi_n, \ph \et_{n+1} = \ga_n$, where $\ga_n \equiv \fr{\xi_n + \et_n}{2}$, if $y_{\ga_n}(t_0) \ge y_{\ga_n}(T)$, otherwise $\xi_{n+1} = \ga_n, \ph \et_{n+1} = \et_n,\ph n=1,2,\dots.$
It is clear from (5.1) and from the definitions of $\{\xi_n\}_{n=1}^{\pinf}$ and  $\{\et_n\}_{n=1}^{\pinf}$ that $\xi_n \le \et_n, \ph \et_n - \xi_n \le \fr{\ga_0}{2^{n-1}}, \ph \xi_n \le \xi_{n+1}, \ph \et_n \ge \et_{n+1}$ and
$$
[y_{\xi_n}(T),y_{\et_n}(T)] \subset [\xi_n,\et_n], \phh n=1,2,\dots .
$$
Therefore $\lim\limits_{n\to\pinf}\xi_n = \lim\limits_{n\to\pinf}\et_n \stackrel {def}{=} \ga_* \ge 0$ and $\bigcap _{n=1}^{\pinf}[y_{\xi_n}(T),y_{\et_n}(T)] = \{\ga_*\}.$  It follows from here that $y_{\ga_*}(t)$ is a nonnegative closed solution of Eq. (1.1). The theorem is proved.

\vsk
By analogy with the proof of Theorem 5.1 it can be proved the following theorem

\vsk
{\bf Theorem 5.2.} {\it Let the conditions (1)--(3) of Theorem 5.1 be satisfied and let \linebreak $d(t) \ge 0, \ph t \in [t_0,T].$ Then Eq. (1.1) has a nonpositive closed solution   on $[t_0,T]$.
}

\phantom{aaaaaaaaaaaaaaaaaaaaaaaaaaaaaaaaaaaaaaaaaaaaaaaaaaaaaaaaaaaaaaaaaaaaaa} $\Box$

{\bf Theorem 5.3.} {\it Let the following conditions be satisfied

\noin
$a(t)< 0$ almost everywhere on $[t_0,T]$,

\noin
$\fr{b^2(t)}{a(t)}$ is integrable on $[t_0,T]$,

\noin
$\il{t_0}{t} \exp\bl\il{t_0}{\ta}\Bigl[c(s) - \fr{b^2(s)}{a(s)}\Bigr]d s\br d(\ta) d \ta \le 0, \phh t\in [t_0,T]$,

\noin
the inequality (2.4) has a solution $\et(t)$ on $[t_0,T]$ such that $\et(t_0) \ge \et(T) > 0$.

\noin
Then Eq. (1.1) has a nonnegative closed solution on $[t_0,T]$.
}

Proof. By virtue of Theorem 3.1 and Remark 3.1 it follows from the conditions of the theorem that for every $\ga \in [0,\et(t_0)]$ Eq. (1.1) has a solution $y_\ga(t)$ on $[t_0,T]$ with $y_\ga(t_0) = \ga$ and $0\le y_\ga(t) \le \et(t), \ph t\in [t_0,T]$.
 Then $0= y_0(t_0) \le y_0(T)$ and $y_{\et(t_0)}(t_0) \ge y_{\et(t_0)}(T)$. Further as in the proof of Theorem 5.1. The theorem is proved.

Using Theorem 3.2 instead of Theorem 3.1 by analogy with the proof of Theorem 5.3 it can be proved the following theorem

{\bf Theorem 5.4.} {\it Let the following conditions be satisfied

\noin
$a(t)< 0$ almost everywhere on $[t_0,T]$,

\noin
$\fr{b^2(t)}{a(t)}$ is integrable on $[t_0,T]$,

\noin
$\il{t_0}{t} \exp\bl\il{t_0}{\ta}\Bigl[c(s) - \fr{b^2(s)}{a(s)}\Bigr]d s\br d(\ta) d \ta \ge 0, \phh t\in [t_0,T]$,

\noin
the inequality (2.8) has a solution $\et(t)$ on $[t_0,T]$ such that $\et(t_0) \le \et(T) < 0$.

\noin
Then Eq. (1.1) has a nonpositive closed solution on $[t_0,T]$.
}

\phantom{aaaaaaaaaaaaaaaaaaaaaaaaaaaaaaaaaaaaaaaaaaaaaaaaaaaaaaaaaaaaaaaaaaaaaa} $\Box$

\vsk

{\bf Theorem 5.5.} {\it Let $y_1(t)$ and $y_2(t)$ be solutions of Eq. (2.1) and Eq. (3.10) respectively on $[t_0,T]$ such that $y_1(t_0) \le y_2(t_0), \ph y_1(t_0) \le y_1(T), \ph y_2(t_0) \ge y_2(T)$ and let the following conditions be satisfied

\noin
 \ph $a(t) < 0$ almost everywhere on $[t_0,T],$ and the function $\fr{b^2(t)}{a(t)}$ is  integrable on $[t_0,T]$,

\noin
 $\ga_1 - y_1(t_0) + \il{t_0}{t}\exp\bl\il{t_0}{\ta}\Bigl[c(s) - \frac{b^2(s)}{a(s)}\Bigr]d s\br\biggl[(a_1(\ta) - a(\ta))y_1^3(\ta) + (b_1(\ta) - b(\ta))y_1^2(\ta) +\linebreak  + (c_1(\ta) - c(\ta)) y_1(\ta) + d_1(\ta) - d(\ta)\biggr] d \ta \ge 0, \ph t\in[t_0,T]$ for some $\ga_1 \in [y_1(t_0),y_2(t_0)]$,

\noin
 $\ga_2 - y_2(t_0) + \il{t_0}{t}\exp\bl\il{t_0}{\ta}\Bigl[c(s) - \frac{b^2(s)}{a(s)}\Bigr]d s\br\biggl[(a_2(\ta) - a(\ta))y_2^3(\ta) + (b_2(\ta) - b(\ta))y_2^2(\ta) +\linebreak  + (c_2(\ta) - c(\ta)) y_2(\ta) + d_2(\ta) - d(\ta)\biggr] d \ta \le 0, \ph t\in [t_0,T]$ for some $\ga_2 \in [\ga_1,y_2(t_0)].$

\noin
Then Eq. (1.1) has a closed solution on $[t_0,T]$
}

Proof. By Theorem 3.3 and Remark 3.1 it follows from the conditions of the theorem that  every solution $y_\ga(t)$ of Eq. (1.1) with $y_\ga(t_0) = \ga \in [y_1(t_0),y_2(t_0)]$ exists on $[t_0,T]$ and
$$
y_1(t_0) \le y_\ga(t) \le y_2(t_0), \phh t \in [t_0,T]. \no (5.4)
$$
Let $\ga_k = y_k(t_0), \ph k=1,2.$ Then it follows from (5.4) and from the conditions  $y_1(t_0) \le y_2(t_0), \ph y_1(t_0) \le y_1(T), \ph y_2(t_0) \ge y_2(T)$ of the theorem that
$$
y_{\ga_1}(t_0) \le y_{\ga_2}(t_0), \ph y_{\ga_1}(t_0) \le y_{\ga_1}(T), \ph \mbox{and} \ph y_{\ga_2}(t_0) \ge y_{\ga_2}(T).
$$
Further as in the proof of Theorem 5.1. The theorem is proved.

\vsk

{\bf Theorem 5.6.} {\it Let $y_1(t)$ and $y_2(t)$ be solutions of Eq. (2.1) and Eq. (3.10) respectively on $[t_0,T]$ such that $y_1(t_0) \le y_2(t_0), \ph y_1(t_0) \ge y_1(T), \ph y_2(t_0) \le y_2(T)$ and let the following conditions be satisfied

\noin
 \ph $a(t) < 0$ almost everywhere on $[t_0,T]$ and the function $\fr{b^2(t)}{a(t)}$ is  integrable on $[t_0,T]$,

\noin
 $\ga_1 - y_1(t_0) + \il{t_0}{t}\exp\bl\il{t_0}{\ta}\Bigl[c(s) - \frac{b^2(s)}{a(s)}\Bigr]d s\br\biggl[(a_1(\ta) - a(\ta))y_1^3(\ta) + (b_1(\ta) - b(\ta))y_1^2(\ta) +\linebreak  + (c_1(\ta) - c(\ta)) y_1(\ta) + d_1(\ta) - d(\ta)\biggr] d \ta \ge 0, \ph t\in[t_0,T]$ for some $\ga_1 \in [y_1(t_0),y_2(t_0)]$,

\noin
 $\ga_2 - y_2(t_0) + \il{t_0}{t}\exp\bl\il{t_0}{\ta}\Bigl[c(s) - \frac{b^2(s)}{a(s)}\Bigr]d s\br\biggl[(a_2(\ta) - a(\ta))y_2^3(\ta) + (b_2(\ta) - b(\ta))y_2^2(\ta) +\linebreak  + (c_2(\ta) - c(\ta)) y_2(\ta) + d_2(\ta) - d(\ta)\biggr] d \ta \le 0, \ph t\in [t_0,T]$ for some $\ga_2 \in [\ga_1,y_2(t_0)].$

\noin
Then Eq. (1.1) has a closed solution on $[t_0,T]$
}

Proof. By virtue of Theorem 3.3 and Remark 3.1 it follows from the conditions of the theorem that every solution $y_\ga(t)$ of Eq. (1.1) with $y_\ga(t_0) = \ga \in [y_1(t_0),y_2(t_0)]$  exists on $[t_0,T]$ and
$$
y_1(t_0) \le y_\ga(t) \le y_2(t), \phh t \in t_0,T]. \no (5.5)
$$
Let $\ga_k = y_k(t_0), \ph k=1,2.$ Then it follows from (5.5) and from the conditions  $y_1(t_0) \le y_2(t_0), \ph y_1(t_0) \ge y_1(T), \ph y_2(t_0) \le y_2(T)$ of the theorem that
$$
y_{\ga_1}(t_0) \le y_{\ga_2}(t_0), \ph y_{\ga_1}(t_0) \ge y_{\ga_1}(T), \ph y_{\ga_2}(t_0) \le y_{\ga_2}(T). \no (5.6)
$$
 Define two sequences $\{\xi_n\}_{n=1}^{\pinf}$ and  $\{\et_n\}_{n=1}^{\pinf}$ by induction on $n$ as follows:

\noin
$\xi_1 = 0, \ph y_{\ga_1}(t_0), \ph \et_1 =  y_{\ga_2}(t_0), \ph \xi_{n+1} =\xi_n, \ph \et_{n+1} = \nu_n$, where $\nu_n \equiv \fr{\xi_n + \et_n}{2}$, if $y_{\nu_n}(t_0) \le y_{\nu_n}(T)$, otherwise $\xi_{n+1} = \nu_n, \ph \et_{n+1} = \et_n,\ph n=1,2,\dots.$
It is clear from (5.6) and from the definitions of $\{\xi_n\}_{n=1}^{\pinf}$ and  $\{\et_n\}_{n=1}^{\pinf}$ that
$$
\xi_n \le \et_n, \ph \et_n - \xi_n \le \fr{\ga_2 - \ga_1}{2^{n-1}}, \ph \xi_n \ge \xi_{n+1}, \ph \et_n \le \et_{n+1}, \ph n=1,2,\dots . \no (5.7)
$$
Therefore $\liml{n\to \pinf}\xi_n = \liml{n \to \pinf}\et_n \stackrel{def}{=} \nu_* \in [\ga_1,\ga_2].$ Show that
$y_{\nu_*}(t)$ is a closed solution of Eq. (1.1) on $[t_0,T]$. It follows from (5.7) that the sequence  $\{y_{\xi_n}(T)\}_{n=1}^{\pinf}$ is non decreasing and upper bounded, the sequence $\{y_{\et_n}(T)\}_{n=1}^{\pinf}$
is non increasing and lower bounded. Moreover, $\liml{n\to \pinf}y_{\xi_n}(T) \le  \liml{n\to \pinf}\xi_n, \ph \liml{n\to \pinf}y_{\et_n}(T) \ge  \liml{n\to \pinf}\et_n.$
It follows from here that $y_{\nu_*}(T) \le y_{\nu_*}(t_0)$ and $y_{\nu_*}(T) \ge y_{\nu_*}(t_0)$. Hence, $y_{\nu_*}(T) = y_{\nu_*}(t_0).$ The theorem is proved.

\vsk

{\bf Corollary 5.1.} {\it Let $y_1(t)$ be a solution of Eq. (2.1) on $[t_0,T]$ such that $y_1(t_0) \le 0$ and let the following conditions be satisfied.

\noin
$a(t) < 0$ almost everywhere on $[t_0,T]$ and the function $\fr{b^2(t)}{a(t)}$ is integrable on $[t_0,T]$,

\noin
$\ga_1 - y_1(t_0) + \il{t_0}{t}\exp\bl\il{t_0}{\ta}\Bigl[c(s) - \frac{b^2(s)}{a(s)}\Bigr]d s\br\biggl[(a(\ta) - a_1(\ta))y_1^3(\ta) - (b_1(\ta) + b(\ta))y_1^2(\ta) +\linebreak  + (c(\ta) - c_1(\ta)) y_1(\ta) - d_1(\ta) - d(\ta)\biggr] d \ta \ge 0, \ph t\in[t_0,T]$ for some $\ga_1 \in [y_1(t_0),y_2(t_0)]$,

\noin
$\ga_2 + y_1(t_0) + \il{t_0}{t}\exp\bl\il{t_0}{\ta}\Bigl[c(s) - \frac{b^2(s)}{a(s)}\Bigr]d s\br\biggl[(a_1(\ta) + a(\ta))y_1^3(\ta) - (b_1(\ta) + b(\ta))y_1^2(\ta) +\linebreak  + (c_1(\ta) + c(\ta)) y_1(\ta) - d_1(\ta) - d(\ta)\biggr] d \ta \le 0, \ph t\in[t_0,T]$ for some $\ga_2 \in [y_1(t_0),-y_1(t_0)]$.

\noin
Then Eq. (1.1) has a closed solution on $[t_0,T]$.
}

Proof. Note that $y_2(t)\equiv -y_1(t)$ is a solution of the equation
$$
y' + a_1(t) y^3 - b_1(t) y^2  + c_1(t) y - d_1(t) = 0, \ph t \ge t_0
$$
on $[t_0,T].$ Then since $Y_1(t_0) \le 0$ by Theorem 5.5 or Theorem 5.6 it follows from the conditions of the corollary that every solution $y(t)$ of Eq. (1.1) with $y(t_0) \in [y_(t_0),-y_1(t_0)]$ exists on $[t_),T]$ and
$$
y_1(t) \le y(t) \le - y_1(t), \phh t \in [t_0,T]. \no (5.8)
$$
Two cases are possible:

\noin
(a) \ph $y_1(t_0) \le y_1(T)$,

\noin
(b) \ph $y_1(t_0) \ge y_1(T)$.

\noin
If the case (a) holds, then $-Y_1(t_0) \ge -y_1(T)$ and by Theorem 5.5 Eq. (1.1) has a closed solution. The case (b) implies $-y_1(t_0) \le -y_1(T)$. In this case existence of a closed solution of Eq. (1.1) on $[t_0,T]$ follows from Theorem 5.6. The corollary is proved.

Using Theorem 3.5 instead of Theorem 3.3 by analogy with the proofs of Theorems 5.5 and 5.6 can be proved the following theorem

{\bf Theorem 5.7} {\it Let $y_1(t)$ and $y_2(t)$ be solutions of Eq. (2.1) and (3.10) respectively on $[t_0,T]$ such that $y_1(t_0) \le y_2(t_0), \ph y_1(t_0) \le y_1(T), \ph y_2(t) \ge y_2(T) \ph (y_1(t_0) \ge y_1(T), \ph y_2(t) \le y_2(T))$ and let the following conditions be satisfied

\noin
$(a_1(t) - a(t)) y_1^3(t) + (b_1(t) - b(t)) y_1^2(t) + (c_1(t) - c(t)) y_1(t) + d_1(t) - d(t) \ge 0, \ph t \in [t_0,T],$

\noin
$(a_2(t) - a(t)) y_2^3(t) + (b_2(t) - b(t)) y_2^2(t) + (c_2(t) - c(t)) y_2(t) + d_2(t) - d(t) \le 0, \ph t \in [t_0,T].$

\noin
Then Eq. (1.1) has a closed solution on $[t_0,T]$.
}

\phantom{aaaaaaaaaaaaaaaaaaaaaaaaaaaaaaaaaaaaaaaaaaaaaaaaaaaaaaaaaaaaaaaaaaaaaa} $\Box$

{\bf Corollary 5.2.} {\it Let $y_1(t)$ be a solution of Eq. (2.1) on $[t_0,T]$ such that $y_1(t_0) \le 0$, and let the following conditions be satisfied

\noin
$(a(t) - a_1(t)) y_1^3(t) - (b_1(t) + b(t)) y_1^2(t) + (c(t) - c_1(t)) y_1(t) - d_1(t) - d(t) \ge 0, \ph t \in [t_0,T],$

\noin
$(a_1(t) + a(t)) y_1^3(t) - (b_1(t) + b(t)) y_1^2(t) + (c_1(t) + c(t)) y_1(t) - d_1(t) - d(t) \le 0, \ph t \in [t_0,T].$

\noin
Then Eq. (1.1) has a closed solution on $[t_0.T]$
}

Proof. By analogy with the proof of Corollary 5.1.

{\bf Example 5.1.} {\it Since $y_1(t)\equiv -1$ is a solution of the equation 
$$
y' + y^3 + \nu(t)y^2 + \nu(t) y +1 = 0, \phh t \ge t_0,
$$
where $\nu(t)$ is a continuous function on $[t_0,\pinf)$, the conditions of Corollary 5.2 are satisfied, if
$$
\sist{a(t) - b(t) + c(t) - d(t) \ge 0,}{-a(t) - b(t) + c(t) - 2 - 2\nu(t) \le 0, \ph t \in [t_0,T].}
$$
Note that we can always chose $\nu(t)$ so large that $-a(t) - b(t) + c(t) - 2 - 2\nu(t) \le 0, \linebreak t~ \in [t_0,T].$ Hence, under the restriction $a(t) - b(t) + c(t) - d(t) \ge 0, \ph t \in [t_0,T]$ Eq. (1.1) has a closed solution on $[t_0,T].$
}

\vsk

\centerline{\bf References}

\vsk

\noin
1. M. Briskin, N. Roytvarf and Y. Yumdin, Center conditions for Abel differential \linebreak \phantom{a} equations. Annals of Mathematics, 172 (2010), 437--483.

\noin
2. M. Briskin and Y. Yumdin, Tangential version of Hilbert 16th problem for the Abel \linebreak \phantom{a}  equation. Moscow Mathematical Journal, vol. 5, No. 1, 2005, 23--53.

\noin
3. V. M. Boyko, Nonlocal symmetry and integrable classes of Abel equation. Proc. of \linebreak \phantom{a} Institute of Mathematics of NAS of Ukraine, 2004, vol. 50, Part 1, 47--51.

\noin
4. Ch.-Sh. Liu, The integrating factor for Riccati and Abel Differential equations. Journal \linebreak \phantom{a}  of  Mathematics Research, Vol. 7, No. 2, 2015, 125--141.

\noin
5. J. Devlin, N. G. Lloid, J. M. Pearson, Cubic systems and Abel equation. J. Differential \linebreak \phantom{a} Equations. 147 (1998(, 435--454.

\noin
6. B. P. Demidovich, Lectures on the mathematical theory of stability. Moskow, \linebreak \phantom{a} "Nauka", 1967.

\noin
7. R. E. Edwards, A formal background to mathematics (Springer-Verlag, New York, \linebreak \phantom{a} Haidelberg, Berlin, 1980).

\noin
8. A. Gassull, J. Libre, Limit cycles for a class of Abel equations. SIAM J. Math. Anal. \linebreak \phantom{a} 21 (1990), 1235--1244.

\noin
9. A. Gassull, R. Prohens, J. Torregrossa, Limit cycles for  rigid cubic systems. J. Math.\linebreak \phantom{a} Anal. Appl.  303 (2005), 391--404.

\noin
10. S. H. Streipert, Abel dynamic equations of the first and second kind. Master \linebreak \phantom{a} Theses,  Student Research and Creative Works, Missuri University of Science and \linebreak \phantom{a} Technology, 2012.

\noin
11. P. J. Torres, Existence of closed solutions for a polynomial first order differential \linebreak \phantom{a}  equation. J. Math. Anal. Appll. 328 (2007) 1108--1116.

\noin
12. A. V.Yurov, A. V. Yaparova, V.A. Yurov, Application of the Abel Equation of 1st  \linebreak \phantom{a} kind to inflation analysis of non-exactly solvable cosmological model. Gravitation and \linebreak \phantom{a} Cosmology, 2014, vol. 20, No. 2, pp. 106--115.

\end{document}